\catcode`\@=11
\newif\if@fewtab\@fewtabtrue
{\count255=\time\divide\count255 by 60
\xdef\hourmin{\number\count255}
\multiply\count255 by-60\advance\count255 by\time
\xdef\hourmin{\hourmin:\ifnum\count255<10 0\fi\the\count255}}
\def\ps@draft{\let\@mkboth\@gobbletwo
    \def\@oddfoot{\hbox to 7 cm{\tiny \versionno
       \hfil}\hskip -7cm\hfil\rm\thepage \hfil {\tiny\draftdate}}
    \def\@oddhead{}
    \def\@evenhead{}\let\@evenfoot\@oddfoot}
\def\draftdate{\number\month/\number\day/\number\year\ \ \ \hourmin }

\global\def\draftcontrol{0}
\def\citen#1{\if@filesw \immediate\write \@auxout {\string\citation{#1}}\fi%
\@tempcntb\m@ne \let\@h@ld\relax \def\@citea{}%
\@for \@citeb:=#1\do {\@ifundefined {b@\@citeb}%
    {\@h@ld\@citea\@tempcntb\m@ne{\bf ?}%
    \@warning {Citation `\@citeb ' on page \thepage \space undefined}}%
    {\@tempcnta\@tempcntb \advance\@tempcnta\@ne
    \setbox\z@\hbox\bgroup\ifcat0\csname b@\@citeb \endcsname \relax
    \egroup \@tempcntb\number\csname b@\@citeb \endcsname \relax
    \else \egroup \@tempcntb\m@ne \fi \ifnum\@tempcnta=\@tempcntb
    \ifx\@h@ld\relax \edef \@h@ld{\@citea\csname b@\@citeb\endcsname}%
    \else \edef\@h@ld{\hbox{--}\penalty\@highpenalty
    \csname b@\@citeb\endcsname}\fi
    \else \@h@ld\@citea\csname b@\@citeb \endcsname \let\@h@ld\relax \fi}%
\def\@citea{,\penalty\@highpenalty\hskip.13em plus.13em minus.13em}}\@h@ld}
\def\@citex[#1]#2{\@cite{\citen{#2}}{#1}}%
\def\@cite#1#2{\leavevmode\unskip\ifnum\lastpenalty=\z@\penalty\@highpenalty\fi%
  \ [{\multiply\@highpenalty 3 #1%
  \if@tempswa,\penalty\@highpenalty\ #2\fi}]}   %
\makeatother 
\catcode`\@=12

\def\be            {\begin{equation}}
\def\bearll        {\begin{array}{ll}}
\def\calc          {\mbox{$\mathcal C$}}
\def\calh          {\mbox{$\mathcal H$}}
\def\caln          {\mbox{$\mathcal N$}}
\def\cft           {conformal field theory}
\def\Cft           {Conformal field theory}
\def\cir           {\,{\circ}\,}

\global\def\draftcontrol{0}
\def\dim           {{\rm dim}}
\def\ee            {\end{equation}}
\def\eear          {\end{array}}
\newcommand\Eev[1] {{{}^{\vee\!}}\!{#1}}
\def\eq            {\,{=}\,}
\newcommand\erf[1] {(\ref{#1})}

\def\fsi           {Fro\-be\-ni\-us\hy Schur indicator}
\def\Hom           {{\rm Hom}}

\def\hy            {$\mbox{-\hspace{-.66 mm}-}$}
\def\id            {\mbox{\sl id}}

\def\iN            {\,{\in}\,} 
\newcommand\labl[1]{\label{#1}\ee \ifnum\draftcontrol=1
                   \mbox{ }\\[-12 mm]\query{#1}\\[5 mm] \fi}
\def\llb           {\mbox{\large[}}
\def\lrb           {\mbox{\large]}}
\def\Nu            {\mbox{$\mathcal V$}}
\def\obj           {{\mathcal O}bj}
\def\one           {{\bf1}}
\def\oti           {\,{\otimes}\,}
\def\Oti           {{\otimes}}

\long\def\query#1{\hskip 0pt{\vadjust{\everypar={}\small\vtop to 0pt{\hbox{}%
     \vskip -13pt\rlap{\hbox to 49.0pc{\hfil{\vtop{\hsize=8pc\tolerance=6000%
     \hfuzz=.5pc\rightskip=0pt plus 3em\noindent#1}}}}\vss}}}}%

\def\twodim        {two-di\-men\-si\-o\-nal}
\def\Vee           {^\vee}

\documentclass[12pt]{article} 
\usepackage{amssymb,amsfonts} 
\usepackage{latexsym}
\begin{document} 
 
\begin{flushright}  {~} \\[-1cm]
{\sf math.QA/0110257}\\[1mm]{\sf PAR-LPTHE 01-46}\\[1mm]
{\sf October 2001} \end{flushright}
 
\begin{center} \vskip 14mm
{\Large\bf A REASON FOR FUSION}\\[3mm]
{\Large\bf RULES TO BE EVEN} \\[15mm]
{\large J\"urgen Fuchs,$\;^1$ \ Ingo Runkel$\;^2$ \ and Christoph 
Schweigert$\;^2$}
\\[8mm]
$^1\;$ Institutionen f\"or fysik, Karlstads universitet~{}\\
Universitetsgatan 1\\ S\,--\,651\,88\, Karlstad\\[5mm] 
$^2\;$ LPTHE, Universit\'e Paris VI~~~{}\\
4 place Jussieu\\ F\,--\,75\,252\, Paris\, Cedex 05
\end{center}
\vskip 18mm
\begin{quote}{\bf Abstract}\\[1mm]
We show that certain tensor product multiplicities in semisimple braided 
sovereign tensor categories must be even. The quantity governing this 
behavior is the Frobe\-ni\-us\hy Schur indicator. The result applies in 
particular to the representation categories of large classes of groups, 
Lie algebras, Hopf algebras and vertex algebras.
\end{quote}
\newpage


The tensor product of $G$-representations is a notion that exists for many
algebraic structures $G$, like groups, Hopf algebras or vertex algebras.
In applications in mathematics and physics, often an important role is played 
by the space of invariants in a tensor product $V_{\!\mu_1}\oti V_{\!\mu_2}
\,\Oti\cdots\Oti\,V_{\!\mu_\ell}$ of irreducible $G$-modules $V_{\!\mu_i}$.

In most cases there is in addition a notion of conjugate representation $V\Vee$.
Irreducible modules which are self-conjugate, i.e.\ satisfy $V\,{\cong}\,V\Vee$,
typically come in two classes -- real (or orthogonal) modules and pseudo-real
(or symplectic, or quaternionic) modules, respectively. They are distinguished 
by the so-called \fsi, which takes value $\nu\eq1$ for real modules and 
$\nu\eq{-}1$ for pseudo-real modules. 

For brevity, we refer to the dimension $\caln_{\!\!\mu_1\mu_2... \mu_\ell}$ 
of the space of invariants in the tensor product 
$V_{\!\mu_1} \oti V_{\!\mu_2}\,\Oti\cdots\Oti\,V_{\!\mu_\ell}$ as a {\em fusion
rule\/}, a terminology borrowed from \twodim\ \cft. When dealing with 
representations of compact Lie groups, it is known that the fusion rule
$\caln_{\!\!\mu_1\mu_2... \mu_\ell}$ is an even integer if all the irreducible
modules $V_{\!\mu_i}$ are self-conjugate and the product 
$\prod_{i=1}^\ell\!\nu_{\mu_i}$ of their \fsi s equals $-1$. In contrast,
when the product of the \fsi s is equal to $+1$, the value of the fusion rule
is not restricted. To mention another example, recent studies of theories of 
unoriented strings and their underlying conformal field theories lead to the 
conjecture \cite{huss2} that 
this relationship also holds for representations of (rational) vertex algebras.

\medskip

In the present letter we show that this property of fusion rules can be 
formulated in a very general setting -- (braided, semisimple) {\em sovereign 
tensor categories\/} -- in which it also finds its natural proof. A sovereign 
tensor category \calc\ is a tensor category with a left and a right duality 
that coincide as functors. We take the tensor category $\calc$ to be strict, 
so that the tensor product is strictly associative.
(By the coherence theorems, this can be assumed without loss of generality.)
For simplicity, we also assume that the morphism sets $\Hom(X,Y)$ of 
$\calc$ are vector spaces over some field $k$.

A right duality assigns to every object $X$ in the tensor category $\calc$
another object $X^\vee$, called the right-dual object, and
gives a family of morphisms
  \be  b_X \in\Hom(\one,X\Oti X\Vee) \qquad{\rm and}\qquad
  d_X \in \Hom(X\Vee\Oti X,\one)  \labl{1d}
($\one$ denotes the tensor unit, satisfying $\one\Oti X\eq X\eq X\Oti\one$) 
such that the equalities
  \be  (\id_X\oti d_X) \cir (b_X\oti\id_X) = \id_X  \qquad{\rm and}\qquad
  (d_X\oti \id_{X\Vee}) \cir (\id_{X\Vee}\oti b_X) = \id_{X\Vee}  \labl{ax}
hold. These data allow to associate to every morphism $f\iN\Hom(X,Y)$ its
right-dual morphism
  \be  f\Vee := (d_Y\oti\id_{X\Vee}) \circ (\id_{Y\Vee}\oti f\oti\id_{X\Vee})
  \circ (\id_{Y\Vee}\oti b_X) \ \in \Hom(Y\Vee\!{,}\,X\Vee) \,.
  \labl{fv}
Similarly, a left duality assigns a left-dual $\Eev X$ and gives morphisms
  \be \tilde b_X \in\Hom(\one,\Eev X\Oti X) \,, \qquad
  \tilde d_X \in \Hom(X\Oti\Eev X,\one)  \end{equation}
that satisfy relations analogous to those in \erf{ax}, and left-dual
morphisms $\Eev f$ are defined analogously as in \erf{fv}.

A {\em sovereign tensor category\/} has both a left and a right duality which 
coincide both on all objects and on all morphisms: $\Eev X\eq X\Vee$ for all
$X\iN\obj(\calc)$, and 
  \be  \Eev f\eq f\Vee \quad \mbox{ for all } \,\, f\iN\Hom(X,Y)  \ee
and all $X,Y\iN\obj(\calc)$. Sovereignty is a strong property; in particular
it allows to define two notions ${\rm tr}_L(f)$ and ${\rm tr}_R(f)$ of a 
trace of an endomorphism $f$, both of which are cyclic. When
applied to the identity morphism, these traces assign a left and a right 
(quantum) dimension $\dim_L(X)\eq{\rm tr}_L(\id_X)$ and
$\dim_R(X)\eq{\rm tr}_R(\id_X)$, respectively, to each object $X$. 
(For our purposes, it is not necessary that the category is spherical, 
i.e.\ that the two traces coincide. For more information and references about 
dualities in tensor categories, see e.g.\ \cite{fuSc16}.) An immediate 
consequence of the definition of sovereignty is the simple, but useful, identity 
  \be \tilde d_Y \circ (g \oti \id_{Y\Vee}) = d_Y \circ (\id_{Y\Vee} \oti g) 
  \labl{eso}
that is valid for all objects $Y$ of $\calc$ and all $g\iN\Hom(\one,Y)$.

\medskip

We now restrict our attention to absolutely simple\,%
\footnote{~In many cases, e.g.\ when \calc\ is abelian and the ground field
$k$ is algebraically closed, the notions of an absolutely simple and of a
simple object are equivalent.}
objects $X$, i.e.\ 
objects for which the ring of endomorphisms coincides with the ground ring, 
$\Hom(X,X)\,{\cong}\,\Hom(\one,\one) \eq k$. In every sovereign tensor 
category there is the following notion of a \fsi\ of an absolutely simple 
self-dual object $X$: Fix an isomorphism $\Phi\iN\Hom(X,X\Vee)$ and consider
  \be  \Nu_X(\Phi) := (d_X\oti\id_{X\Vee}) \circ
  (\id_{X\Vee}\oti\Phi^{-1}\oti\Phi) \circ (\id_{X\Vee}\oti\tilde b_X)
  \;\in \Hom(X\Vee\!{,}\,X\Vee) \,. \labl{nu}
This morphism does in fact not depend on $\Phi$, and it must be a multiple of 
the identity,
  \be  \Nu_X(\Phi) = \nu_X\, \id_{X^\vee} \,.  \ee
The scalar $\nu_X$ is called the {\em\fsi\/} of $X$. If all left and right 
dimensions in the category are invertible, then $\nu_X$ can only take the
values $\pm1$. In the definition of $\Nu_X$ one may exchange the role of 
left and right dualities, but by the sovereignty of \calc\ both definitions 
are equivalent. The following identity, valid for all morphisms 
$\Phi\iN\Hom(X,X\Vee)$, follows immediately from the definitions:
  \be  \tilde d_X \circ (\id_X\otimes\Phi)
  = \nu_X \, d_X \circ (\Phi\otimes \id_X) \,.  \ee

\medskip

Having explained the setting, we can formulate our main result:

\medskip \noindent
{\bf Theorem:}\,
{\sl Let $\calc$ be a semisimple braided sovereign tensor category, 
$\one\iN\obj(\calc)$ the tensor unit, and $X_i\iN\obj(\calc)$ {\rm(}$i\iN
\{1,2,...\,,\ell\}${\rm)} self-dual absolutely simple objects with \fsi s 
$\nu_i\iN\{\pm1\}$. Then the morphism space
  \be  \calh := \Hom(\one,X_1\oti X_2\,\Oti\cdots\Oti\,X_\ell) \ee
can be endowed with a non-degenerate bilinear pairing.
\\[.1em]
This pairing is symmetric if the product $\nu\,{:=}\,\prod_{i=1}^\ell\!\nu^{}_i$
is $+1$, and antisymmetric if $\nu\eq{-}1$. As a consequence, for
$\nu\eq{-}1$ the vector space \calh\ is even-dimensional.}

\medskip

As a matter of fact, we shall prove a slightly more general statement.
To formulate it, it is convenient to define the \fsi\ of a non-selfdual 
absolutely simple object $X$ to be zero, $\nu_X\eq0$ for $X\,{\ncong}\,X^\vee$.

\medskip \noindent
{\bf Theorem:}\,
{\sl Let \calc\ be a semisimple
braided sovereign tensor category, $\one\iN\obj(\calc)$ the tensor unit, and
$X_i\iN\obj(\calc)$ {\rm(}$i\iN\{1,2,...\,,\ell\}${\rm)} absolutely simple
objects with \fsi s $\nu_i\iN\{0,\pm1\}$. Suppose there exists a
permutation $\pi\iN S_\ell$ such that $X_{\pi(i)}\,{\cong}\,X_i^\vee$
for $i\eq1,...\,,\ell$. Then the morphism space
  \be  \calh := \Hom(\one,X_1\oti X_2\,\Oti\cdots\Oti\,X_\ell)  \ee
can be endowed with a non-degenerate bilinear pairing.
\\[.1em]
This pairing is symmetric if the product $\nu\,{:=}\prod_{i=1}^\ell(1+\nu_i
-\nu_i^2)$ is $+1$, and antisymmetric if $\nu\eq{-}1$. In particular, for
$\nu\eq{-}1$ the vector space \calh\ is even-dimensional.}

\medskip

\noindent
{\bf Proof}: \\
1.\ The presence of the permutation $\pi$ allows us to find another permutation 
$\sigma\iN S_\ell$ of order two satisfying $X_{\sigma(i)}\,{\cong}\,X_i^\vee$
as well. \\
To see this, we decompose $\pi$ in disjoint cycles $(i_1,...\,,i_n)$. On each 
such cycle we have $X_{i_k}\,{\cong}\,X_{i_{k+1}}^\vee$\linebreak[0]${\cong}\,
X_{i_{k+2}}$ 
(we cyclicly identify $i_{n+1}\eq i_1$ and $i_{n+2}\eq i_2$). If the length $n$
of a cycle is odd, then $X_{i_k}\,{\cong}\,X_{i_{k+1}}$ and hence all 
objects on the cycle are self-dual. On such a cycle we let $\sigma$ act as 
the identity. If the order $n$ is even, we take $\sigma(i_k)\,{:=}\,i_{k+1}$ 
for even $k$ and $\sigma(i_k)\,{:=}\,i_{k-1}$ for odd $k$.
Note that $\sigma(i)\,{\ne}\, i$ does not rule out that $X_i$ is self-dual. 
\\[.2em]
2.\ Next we fix isomorphisms $f_i \iN \Hom(X_i, X_{\sigma(i)}^\vee)$: If 
$\sigma(i)\eq i$, we pick an arbitrary isomorphism $f_i\iN \Hom(X_i, X_i^\vee)$.
For $\sigma$-orbits of length 2, we choose an arbitrary isomorphism $f_i\iN 
\Hom(X_i, X_{\sigma(i)}^\vee)$ for one element $i$ on the orbit, and then define
$f_{\sigma(i)} \iN \Hom(X_{\sigma(i)}, X_i^\vee)$ as
  \be  f_{\sigma(i)} := (id_{X_i^\vee} \oti d_{X_{\sigma(i)}})
  \circ ( \id_{X_i^\vee} \oti f_i \oti \id_{X_{\sigma(i)}})
  \circ (\tilde b_{X_i} \oti \id_{X_{\sigma(i)}}) \,.  \labl{isodef}
(Notice that $f_{\sigma(i)}$ is not the (left- and right-) dual morphism
of $f_i$.) This definition does not depend on the choice of $i$ on the
orbit. Indeed, formula \erf{isodef} implies that
  \be  f_i = (id_{X_{\sigma(i)}^\vee} \oti d_{X_i}) \circ
  ( \id_{X_{\sigma(i)}^\vee} \oti f_{\sigma(i)} \oti \id_{X_i}) \circ
  (\tilde b_{X_{\sigma(i)}} \oti \id_{X_i}) \,,  \ee
which is shown by applying \erf{eso} to the morphism $g\,{:=}\,(\id_{X_i\Vee}
\Oti f_i)\cir\tilde b_{X_i}\iN\Hom(\one,X_i\Vee\Oti X\Vee_{\sigma(i)})$.
With this choice of the isomorphisms $f_i$, we have
  \be  d_{X_{\sigma(i)}} \circ (f_i\oti\id_{X_{\sigma(i)}})
  = p_i \; \tilde d_{X_i} \circ (\id_{X_i}\oti f_{\sigma(i)}) \,, \labl{f-prop}
where $p_i\eq\nu_i$ if $\sigma(i)\eq i$ and $p_i\eq1$ if $\sigma(i)\,{\ne}\,i$.
\\[.2em] 
3.\ We now construct a bilinear pairing $\langle\,\cdot\,,\cdot\,\rangle{:}\,\
\calh\,{\times}\,\calh\:{\to}\;k\,{\equiv}\,\Hom(\one,\one)$ 
out of the following two ingredients: First, the duality morphism
  \be  d_{X_{\sigma(1)}\otimes X_{\sigma(2)}\otimes\cdots\otimes X_{\sigma(\ell)}}
  \in \Hom(X_{\sigma(1)}\Vee\Oti X_{\sigma(2)}\Vee\Oti\cdots\Oti X_{\sigma(\ell)
  }\Vee
  \oti X_{\sigma(\ell)}\Oti X_{\sigma({\ell-1})}\Oti\cdots\Oti X_{\sigma(1)},
  \one) \,,  \ee
which may be obtained recursively by
  \be  d_{X_{\sigma(1)}\otimes X_{\sigma(2)}\otimes\cdots\otimes X_{\sigma(\ell)}}
  \eq d_{X_{\sigma(1)}}\cir \llb
  \id_{X_{\sigma(1)}\Vee}\oti d_{X_{\sigma(2)}\otimes X_{\sigma(3)}\otimes
  \cdots\otimes X_{\sigma(\ell)}}\oti\id_{X_{\sigma(1)}} \lrb  \ee
from the duality morphisms of the simple objects $X_i$. Second, any isomorphism
  \be  c^{}_{X_1\otimes X_2\otimes\cdots\otimes X_\ell} \in
  \Hom(X_1\Oti X_2\Oti\cdots\Oti X_\ell,X_{\sigma(\ell)}
  \Oti X_{\sigma(\ell-1)}\Oti\cdots\Oti X_{\sigma(1)})  \ee
that is a combination of braidings. (Any of the various possible combinations 
of braidings may be chosen; the argument does not depend on this choice.)
Then the pairing $\langle\,\cdot\,,\cdot\,\rangle$ on \calh\ is given by
  \be  \langle\phi\,,\phi'\rangle
  := d_{X_{\sigma(1)}\otimes X_{\sigma(2)}\otimes\cdots\otimes X_{\sigma(\ell)}}
  \circ \llb f_1\Oti f_2\,\Oti\cdots\Oti f_\ell
  \oti c^{}_{X_1\otimes X_2\otimes\cdots\otimes X_\ell} \lrb
  \circ \llb \phi\oti\phi' \lrb \,.  \labl{<>}
\vskip.1em\noindent
4.\ The pairing \erf{<>} is non-degenerate. \\
Since both $c^{}_{X_1\otimes X_2\otimes\cdots\otimes X_\ell}$ and
$f_1\Oti f_2\,\Oti\cdots\Oti f_{\ell}$ are ismorphisms, this follows 
from the fact that in every semisimple sovereign category for any object $W$
the pairing
  \be  B: \quad \begin{array}{rl}
  \Hom(\one,W\Vee) \oti \Hom(\one,W) \! &\to\; k \\[.1em]
  \phi_1 \otimes \phi_2 \! &\mapsto\, d_W \cir (\phi_1\oti \phi_2)
  = \phi_1\Vee\cir\phi_2  \end{array} \ee
is non-degenerate. Namely, as $\calc$ is semisimple, the object $W$ can
be written as a direct sum $W\,{\cong}\, \bigoplus_i W_i$ of simple objects.
Thus when $\phi_2\iN\Hom(\one,W)$ is non-vanishing, then at least one component 
corresponding to some $W_i\,{\cong}\,\one$ is non-vanishing, and there exists a
$\phi_1$ with a matching component  such that $B(\phi_1,\phi_2)$ is non-zero.
\\[.2em]
5.\ We finally analyze the symmetry properties of the pairing. \\
Using property \erf{f-prop} of the morphisms $f_i$ as well as
functoriality of the braiding, we obtain
  \be  \langle\phi\,,\phi'\rangle
  = \nu\, \tilde d_{X_{\sigma(\ell)}\otimes X_{{\sigma(\ell-1)}}\otimes\cdots
  \otimes X_{\sigma(1)}}\circ\llb c^{}_{X_1\otimes X_2\otimes\cdots\otimes X_\ell}
  \oti f_1\Oti f_2\,\Oti\cdots\Oti f_\ell \lrb \circ \llb\phi\oti\phi'\lrb
  \labl{z2}
with $\nu\eq\prod_{i=1}^\ell\! p_i$. We wish to rewrite $\nu$ in terms of 
\fsi s. For $\sigma(i)\eq i$ we have $p_i\eq\nu_i\iN\{\pm1\}$ and thus 
$p_i\eq1+\nu_i-\nu_i^2$. For $\sigma(i)\,{\ne}\,i$ we have $p_i\eq p_{\sigma(i)}
\eq1$, and in all three possible cases -- that is, $X_i\,{\cong}\,X_i^\vee$
with $\nu\eq1$ or with $\nu\eq{-}1$, or $X_i\,{\ncong}\,X_i^\vee$ -- one has
$(1+\nu_i-\nu_i^2)(1+\nu_{\sigma(i)}-\nu_{\sigma(i)}^2)\eq1\eq p_i\,p_{\sigma(i)}$.
Altogether we obtain $\nu\eq\prod_{i=1}^\ell\!(1+\nu_i-\nu_i^2)$.
Applying the identity \erf{eso} to the expression \erf{z2}, we thus find
  \be  \langle\phi\,,\phi'\rangle = \nu\; \langle\phi',\phi\rangle \,,  \ee
which proves the assertion about the (anti-)symmetry of the pairing.
\hfill $\Box$ 

\bigskip

\noindent
{\bf Remark}: \\
While sovereignty of \calc\ is a crucial ingredient of the proof, the r\^ole 
of semisimplicity and of the braiding is limited. 
\\[.15em]
1.\ Instead of semisimplicity, a substantially weaker property is sufficient: A 
weak form of {\em dominance\/} \cite{TUra}, namely the existence of a family 
$I$ of absolutely simple objects such that the identity morphism $\id_W$ of 
any object $W$ can be decomposed in a finite 
sum $f\eq\sum_r g_r\cir h_r$ with $h_r\iN\Hom(W,i)$ and $g_r\iN\Hom(i,W)$ for
some $i\eq i(r)\iN I$.
\\[.15em]
2.\ Concerning the braiding, all that is needed in the proof is the existence
of an isomorphism
$\tilde c\iN \Hom(X_1\Oti X_2\Oti\cdots\Oti X_\ell,X_{\sigma(\ell)}
\Oti X_{\sigma(\ell-1)}\Oti\cdots\Oti X_{\sigma(1)})$ such that 
  \be  (f_{\sigma(\ell)}\oti f_{\sigma(\ell-1)}\,\Oti\cdots\Oti\,f_{\sigma(1)}) 
  \circ \tilde c
  = \tilde c\Vee {\circ}\, (f_1\oti f_2\,\Oti\cdots\Oti\,f_{\ell}) \,.   \ee
In a general sovereign tensor category without braiding, an isomorphism with
this property need not exist.
In fact, without any further structure there is even no reason why the two 
tensor products in question should be isomorphic at all. However, in every 
sovereign tensor category there is a class of tensor products to which the 
theorem applies immediately: Those for which $X_j\Vee\cong X_{\ell-j+1}$ 
for all $j\eq1,2,...\,,\ell$. In these cases, one may simply take 
$\tilde c\eq\id\,$ and $\sigma(j)\eq\ell-j+1$.

\bigskip

Our theorem applies in particular to the category of finite-dimensional
representations of a finite group or, more generally, of a finite-dimensional
compact Lie group. In this case the statement is well-known at least to experts.  
Apart from these classical applications, our result applies also to such
modular tensor categories \cite{TUra} in which the ground ring 
$\Hom(\one,\one)$ is a field. Examples of modular tensor categories are the 
representation category 
of the quantum double of a finite group, or the (truncated) representation 
category of the deformed enveloping algebra of a simple Lie algebra (a quantum 
group) with the deformation parameter being a root of unity.

Many more examples of modular tensor categories are supplied by rational 
conformal field theories \cite{mose3,muge7}. In this context, the \fsi\ also 
shows up when one studies the conformal field theory on a Klein bottle (see 
\cite{bant11,fffs3}). Correlators on the Klein bottle appear in the 
so-called orientifold projection (see e.g.\ \cite{sasT2}). 
As observed in \cite{huss2}, consistency of this projection 
requires certain coupling spaces to be even-dimensional. This observation 
amounts to the statement of our theorem for representations of rational
vertex algebras, which was our original motivation to study this issue. 


\vskip2.8em \noindent{\bf Acknowledgment}.
It is a pleasure to thank M.\ Duflo, G.\ Mason, J.\ Michel, A.N.\ Schellekens,
J.-P.\ Serre and W.\ Soergel for helpful correspondence.
IR is supported by the EU grant HPMF-CT-2000-00747.
\bigskip

 \newcommand\wb{\,\linebreak[0]} \def\wB {$\,$\wb}
 \newcommand\Bi[1]    {\bibitem{#1}}
 \newcommand\J[5]     {{\sl #5}, {#1} {#2} ({#3}) {#4} }
 \newcommand\K[6]     {{\sl #6}, {#1} {#2} ({#3}) {#4}}
 \newcommand\Prep[2]  {{\sl #2}, pre\-print {#1}}
 \newcommand\BOOK[4]  {{\sl #1\/} ({#2}, {#3} {#4})}
 \def\jf    {J.\ Fuchs}
 \def\dim   {dimension} 
 \def\comp  {Com\-mun.\wb Math.\wb Phys.}
 \def\cpma  {Com\-pos.\wb Math.}
 \def\foph  {Fortschritte\wB d.\wb Phys.}
 \def\phlb  {Phys.\wb Lett.\ B}
 \def\nupb  {Nucl.\wb Phys.\ B}
 \def\NY     {{New York}}


\end{document}